\documentclass[12pt]{article}
\usepackage{amsmath, amsthm, amsfonts, amscd, xypic}

\def\FF{\mathbb{F}}
\def\NN{\mathbb{N}}
\def\QQ{\mathbb{Q}}
\def\ZZ{\mathbb{Z}}

\def\Fpbar{\overline{\FF_p}}
\def\map#1#2#3{#1\!:\!#2 \to #3}

\def\pibar{\tilde{\pi}}

\newtheorem{prop}{Proposition}
\newtheorem{lemma}{Lemma}
\newtheorem{theorem}{Theorem}

\begin{document}

\title{Power Series and $p$-adic Algebraic Closures}
\author{Kiran S. Kedlaya}
\date{\today}

\maketitle

\begin{abstract}
We describe a presentation of
the algebraic closure of the ring of Witt vectors of an
algebraically closed field of characteristic $p>0$. The construction
uses ``generalized power series in $p$'' as constructed by Poonen, based
on an example of Lampert.
\end{abstract}

\section{Introduction}

In a previous paper \cite{me}, the author gave an explicit description of
the algebraic closure of the power series field over a field of
characteristic $p>0$, in terms of certain ``generalized power series''.
The purpose of the present paper is to extend this work to
mixed characteristic. Specifically, we give an analogous
description of the algebraic closure of the
Witt ring $W(K)$ of an algebraically closed field $K$ of characteristic $p$.
In place of generalized power series, we use ``generalized $p$-adic series''
as introduced by Poonen \cite{bjorn}.
(Note: here and throughout, the ``algebraic closure'' of a domain refers to
the integral closure of the domain in the algebraic closure of its 
function field.)

Our approach is to relate the algebraic closure of $W(K)$ 
to the algebraic closure of $K[[t]]$ via the Witt ring of the latter. 
We first exhibit a surjection of $W(\overline{K[[t]]})$ 
onto $\overline{W(K)}^\wedge$ (the wedge denotes $p$-adic
closure) in which $t$ maps to $p$. (The wedge denotes
$p$-adic completion.) Using the results of 
\cite{me}, we then give explicit presentations of 
$W(\overline{K[[t]]})$ and of $\overline{W(K)}^\wedge$.

Our approach makes use of two notions which occur in \cite{me} 
but might otherwise be unfamiliar to the reader: generalized power 
series, described in Section~\ref{sec:def}, and twist-recurrent 
sequences (or linearized recurrent sequences), described in 
Section~\ref{sec:twi}.


\section{Generalized power series} \label{sec:def}

In this section, we describe the construction of generalized power
series and their $p$-adic analogues; before doing so, however, let us
briefly review some of the history of these constructions, following
the account given in \cite{bjorn}.

Generalized power series were first introduced by Hahn \cite{hahn},
and were studied in terms of valuations by Krull \cite{krull}.
(Nowadays rings of generalized power series are sometimes known as
Mal'cev-Neumann rings, as these two authors independently extended
the construction to nonabelian value groups. We will not need that
generalization in this paper.)
It was noted by Kaplansky \cite{kap} that his ``maximal immediate
extension'' of a valued field turns out to be a
ring of generalized power series 
in the equal
characteristic case.
The mixed-characteristic analogues of generalized power 
series were introduced by Poonen \cite{bjorn},
who showed that the ring they form realizes Kaplansky's maximal immediate
extension in the mixed characteristic case.
Poonen's construction was motivated by an example of Lampert \cite{lamp}.

Let us now proceed to the  constructions.
For any ring $R$ and any totally ordered abelian group $G$ (whose
identity element we call 0), we define
the set $R((t^G))$ of formal power series over $R$ with value group
$G$ as the set of functions $\map{x}{G}{R}$ for which the 
set of $i \in G$ with $x_i \neq 0$ is \emph{well-ordered}, that is,
contains no
infinite decreasing subsequence; this set is called the
\emph{support} of $x$.
We will often write the element $x$ as $\sum_i 
x_i t^i$, where $t$ is the dummy variable specified in the notation 
$R((t^G))$. 
With this notation, $R((t^G))$ begs to be given the structure of a ring with 
the operations
\begin{align*}
\sum_i x_it^i + \sum_j y_j t^j &= \sum_k (x_k+y_k) t^k \\
\sum_i x_it^i \times \sum_j y_j t^j &= \sum_k \left(\sum_{i+j=k} x_i 
y_j\right) t^k.
\end{align*}
Fortunately, these definitions make sense, the former because the union 
of two well-ordered sets is well-ordered, the latter because the set 
of sums of elements of two well-ordered sets is well-ordered and 
because any such sum can be so expressed in finitely many ways. (These 
are all easy consequences of the no-infinite-decreasing-subsequence 
definition.)

The subring $R[[t^G]]$ of $R((t^G))$ comprises those $x$ whose 
support consists entirely of nonnegative elements of $G$. Both of 
these rings carry a natural $t$-adic valuation $v_t$ with values in $G$.

\begin{prop} \label{divis}
Assume $G$ is divisible. If $R$ is an algebraically closed field,
then so is $R((t^G))$.
\end{prop}
We give a proof using a transfinite version of Newton's algorithm;
another proof can be found in
\cite{rib}, and a proof using the theory of immediate extensions is sketched in
\cite{bjorn}.
\begin{proof}
Associated to the polynomial
$P(x) = a_0 x^n + a_1 x^{n-1} + \cdots + a_n$ over $R((t^G))$ is its
\emph{Newton polygon}, the lower boundary of the convex hull
of the points $(i, v_t(a_i))$ for $i=0, \dots, n$. The integers $m$ such that
$(m, v_t(a_m))$ are vertices of the Newton polygon are called the \emph{breakpoints},
and the ratios $(v_t(a_{m_1}) - v_t(a_{m_2}))/(m_1 - m_2)$, where $m_1, m_2$ are
adjacent breakpoints, are called the \emph{slopes} of the Newton polygon (they exist
because $G$ is divisible).
We will eventually see
that the valuations of the roots of $P$ are the slopes of the constituent segments
of the Newton polygon, and keeping this in mind will clarify the motivation of the following argument.

Let $\omega_G$ be the smallest ordinal with cardinality greater than that of $G$.
We show that for any monic polynomial $P(x)$ over $R((t^G))$,
there exists a map $f: \omega_G \to G \cup \{\infty\}$
and a map $g: \omega_G \to R$ such that $r_\omega =
\sum_{\alpha < \omega} g(\alpha) t^{f(\alpha)}$ has the following properties.
(In the definition of $r_\omega$, we formally take $t^\infty = 0$.)
\begin{enumerate}
\item
If $\omega_1 < \omega_2$, then $f(\omega_1) \leq f(\omega_2)$, with equality only
if both are equal to $\infty$.
\item
If $f(\omega) = \infty$, then $r_\omega$
is a root of $P$.
\item
If $f(\omega) < \infty$, then the polynomial $P(x-r)$ has largest Newton slope strictly
greater than $f(\alpha)$.
\end{enumerate}
Since there is no injective map from $\omega_G$ to $G$,
there exists some $\omega < \omega_G$ such that $f(\omega) = \infty$, and the corresponding
$r_\omega$ will be a root of $P$, proving the theorem.
In other words, we will carry out a ``transfinite Newton's algorithm'' and show that
it converges to a root of $P$.

We prove the claim by transfinite induction. The induction step is self-evident for limit
ordinals, so we may need only worry about non-limit ordinals as well as the base case. Both
of these are subsumed in the following fact: if $P(0) \neq 0$,
there exist $r \in R$ and $s \in G$ such that the largest slope of the Newton polygon
of $P(x - rt^s)$ is strictly greater than the largest slope of the Newton polygon of $P(x)$.

To prove the latter claim, write $P(x) = a_0x^n + a_1 x^{n-1} + \cdots + a_n$
with $a_0 = 0$.
Let $m$ be the largest breakpoint (we do
not regard $n$ itself as a breakpoint) and $s$ the slope of the segment of the Newton
polygon between $m$ and $n$. Since $K$ is algebraically closed, the polynomial
\[
A(x) = a_m t^{-v_t(a_m)} x^{n-m} + a_{m+1} t^{-v_t(a_m)-s} x^{n-m-1} + \cdots +
a_n t^{-v_t(a_m)-s(n-m)}
\]
factors completely over $K$. Let $r$ be any root of this polynomial, let $q$ be
the multiplicity of the root $r$, and put
$Q(x) = P(x-rt^s)$. If we put $Q(x) = \sum b_i x^{n-i}$, then
\[
b_i = \sum_j a_{i-j} \binom{n-i+j}{j} r^j t^{sj}.
\]
From this description we can read off the Newton slopes of $Q(x)$.
First note that for any breakpoint
$i \leq m$, the sum defining $b_i$ consists of
$a_i$ plus various terms of valuation strictly larger than that of
$a_i$.
For $i > m$, we have $v_t(b_i) \geq v_t(b_m) + s(i-m)$, and
the coefficient of $t^{v_t(b_m) + s(i-m)}$ in $b_i$ is equal to the
coefficient of $x^{i-m}$ in $A(x-r)$. In particular, we have
$v_t(b_{n-q}) = v_t(b_m) + s(n-m-q)$ and $v_t(b_i) > v_t(b_m) + s(i-m)$
for $i > n-q$. In short, the Newton slopes of $Q(x)$ less than $r$ are
the same as those of $P$, the slope $r$ occurs with multiplicity $n-m-q$,
and the remaining slopes are greater than $r$. This completes the proof of the
claim and hence of the induction step.
\end{proof}

If the ring $R$ carries a valuation $v$, we define
the subring $R_v((t^G))$ of $R((t^G))$ as the set of $x$ such that
for any $r$, the set of $i \in G$ such that 
$v(x(i)) \leq r$ is well-ordered.  This ring 
is canonically isomorphic to the inverse limit of 
$R/I_r((t^G))$ over all $r$, where $I_r = \{t \in R: v(t) \geq r\}$.
In particular, if $v$ 
is a discrete valuation, $R_v((t^G))$ is $I_1$-adically complete.

We now proceed to the promised  mixed-characteristic analogue.
Given a prime $p$, a $p$-adically complete
ring $R$, and a totally ordered abelian group $G$ equipped with a fixed 
embedding of $\ZZ$ in $G$,
the ring $R((p^G))$ of $p$-adic series over $R$ with value group $G$
is the quotient
of $R((t^G))$ by the ideal consisting of those $x$ for which
$\sum_{n=0}^\infty x_{n+i} p^n = 0$ for all $i \in G$;
the ring $R[[p^G]]$ is defined as the analogous quotient of $R[[t^G]]$.
In case $R/pR$ is perfect, we can choose canonical 
representatives of elements of $R((p^G))$ in $R[[t^G]]$, namely those 
whose coefficients are all Teichm\"uller elements of $R$.

Unsurprisingly, we have the following analogue of
Proposition~\ref{divis}. The proof of Proposition~\ref{divis}
given in \cite{bjorn} applies in mixed characteristic as well;
alternatively, one can adapt the proof given above simply by
replacing $t$ by $p$ everywhere it appears.
\begin{prop} \label{divis2}
Assume $G$ is divisible. If $R$ is an algebraically closed field,
then so is $W(R)((p^G))$.
\end{prop}

\section{Twist-recurrent sequences} \label{sec:twi}

Throughout this section, let $R$ be a $p$-adically complete and separated
ring whose residue ring is an algebraically closed field,
and $\map{\sigma}{R}{R}$ an automorphism of $R$ inducing the Frobenius
automorphism $x \mapsto x^p$
on the residue field. 
(We sometimes write $F$, acting as an operator on the left, instead of
$\sigma$, acting as a superscript on the right.)
Let $R_0$ be the fixed ring of $\sigma$, which
is complete with residue field $\ZZ_p$, and let $\pi$ be a generator of the
maximal ideal of $R_0$, which then also generates the maximal ideal of $R$.

The main case of interest is when $R$ is the ring of Witt vectors
of an algebraically closed field of characteristic $p$,
or the quotient of said ring by the ideal
generated by a power of $p$.
The bulk of this section follows \cite[proof of 
Lemma~4]{me}, \textit{mutatis mutandis}; significant
departures are noted where they occur.

A \emph{twist-recurrence relation} is an equation of the form
\begin{equation} \label{twist}
d_0 c_n + d_1 c_{n+1}^\sigma + \cdots + c_{n+k}^{\sigma^k} = 0 \qquad
\forall n \geq 0,
\end{equation}
where $d_0, \dots, d_k \in R$ are not all zero, $d_0$ is
not divisible by $\pi$ and $c_0, c_1, \dots$ 
is an infinite sequence of elements of $R$; a \emph{twist-recurrent 
sequence} is any sequence of the form $\{c_n\}$ for suitable $k$ and 
$d_0, \dots, d_k$.

Before proceeding further, we note that in \cite{me}, it was not required that $d_0$ be
nonzero, and so a twist-recurrent sequence with finitely many terms prepended
remained twist-recurrent. In this more general setting, allowing $d_0$ to be
divisible by $\pi$ would cause technical complications not relevant to our main task.

The basic properties of twist-recurrent sequences superficially
resemble those of
ordinary linear recurrent sequences (those satisfying the same
defining equation with the $\sigma$-action removed).
For example, it is evident that
(assuming $d_k \neq 0$) the set of such sequences forms an $R$-module
with scalar multiplication defined by
\[
\lambda \{c_n\} = \{ \lambda^{\sigma^{-n}} c_n \}
\]
(but not under the usual scalar multiplication $\lambda \{c_n\} = \{ 
\lambda c_n\}$).

To continue the analogy, we note that twist-recurrent sequences can be 
studied by exhibiting a canonical basis for the set of sequences 
satisfying a given relation. Of course, the exact nature of these 
basis sequences is somewhat different than in the classical case:
they are the constant sequences $c_n = z$ for $z$ a root of the equation
\[
P(F)(z) = (d_0 + d_1 F + \cdots + d_n F^n)z = 0.
\]
The following lemma (which may be thought of as
Hensel's Lemma for polynomials in 
the operator $F$)
shows that these sequences indeed furnish a 
complete basis for the set of solutions of a twist-recurrence relation.
\begin{lemma}
Let $R_0$ be the fixed ring of $\sigma$.
Then given $d_0, \dots, d_{n-1} \in R$ with $d_0 \not\equiv 0 \pmod{\pi}$, 
the set of solutions of the equation
\[
P(F)z = (F^n + d_{n-1} F^{n-1} + \cdots + d_0)z = 0
\]
in $R$ is a free $R_0$-module of rank $n$.
\end{lemma}
\begin{proof}
Let $z^{(1)}_1, \dots, z^{(1)}_n \in R$ be elements which reduce to
a basis over $\FF_p$ of
the $\FF_p$-vector space solutions of $P(F)z=0$ modulo $\pi$.
For each $k \in \NN$, we wish to exhibit $z_1^{(k)}, \dots, z_n^{(k)}$ which
reduce to a basis over $R_0/\pi^n R_0$ of the set of solutions of
$P(F)z=0$ in $R/\pi^n R$. We construct these by an inductive argument:
given the $z_i^{(k)}$, note that
\[
P(F)(z_i^{(k)} + s\pi^k) \equiv P(F)z_i^{(k)} + \pi^k(s^{p^n} + d_{n-1} s^{p^{n-1}} + 
\cdots + d_0 s) \pmod{\pi^{k+1}}.
\]
Since $R$ has algebraically closed residue field, there exists $s \in R$
such that $P(F)(z_i^{(k)} + s\pi^k \equiv 0 \pmod{\pi^{k+1}}$;
we set $z_i^{(k+1)} = z_i^{(k)} + s\pi^k$.
The facts that the $z_i^{(k)}$ are independent and that they span the
solution space of $P(F)z=0$ are straightforward.
\end{proof}
Note that the lemma can (and in fact always will) fail if $d_0$ is 
divisible by $\pi$; for example, the equation $(F+\pi)z = 0$ has no 
nonzero solutions. However, in case $\pi$ is not nilpotent,
a full set of solutions will exist if we 
adjoin $\pi^{1/m}$ and a solution of 
$(F+\pi^{1/m})z = 0$ for some $m$, and the following lemma (the analogue 
of \cite[Corollary~5]{me} in this context) can be proved 
in this fashion without the restriction on $d_0$.

\begin{lemma} \label{twistcomb}
For $d_{0}, \dots, d_{k-1}, d'_{0}, \dots, d'_{l-1} \in R$
with $d_0, d'_0 \not\equiv 0 \pmod{\pi}$, consider
the twist-recurrence relations
\begin{align*}
c_{n+k}^{\sigma^k} + d_{k-1} c_{n+k-1}^{\sigma^{k-1}} + \cdots + d_{0}
c_{n} &= 0 \\
(c'_{n+l})^{\sigma^l} + d'_{l-1} (c'_{n+l-1})^{\sigma^{l-1}} + \cdots + d'_{0}
c'_{n} &= 0.
\end{align*}
Then the following are true:
\begin{enumerate}
\item
There is a twist-recurrence relation whose solutions include all sequences of
the form $(c_{n} + c'_{n})$, where $(c_{n})$ satisfies the first
relation and $(c'_{n})$ satisfies the second.
\item
There is a twist-recurrence relation whose solutions include all sequences of
the form $(c_{n}c'_{n})$, where $(c_{n})$ satisfies the first
relation and $(c'_{n})$ satisfies the second.
\end{enumerate}
\end{lemma}
\begin{proof}
By the previous lemma, the sequences $(c_n)$ and $(c'_n)$ satisfying 
the given relations have the form
\[
c_n = \sum_i \lambda_i^{\sigma^{-n}} z_i, \quad
c'_n = \sum_j \lambda_j^{\sigma^{-n}} z'_j,
\]
where the $z_i$ and $z'_j$ depend only on the $d_i$ and $d'_j$.
Clearly $(c_n + c'_n)$ and $(c_n c'_n)$ are also of this form, with 
the set of $z_i$ being the union of the sets of $z_i$ and $z'_j$ in 
the former case, and the set of products $z_i z'_j$ in the latter case.

Thus it suffices to show that given a set of $z_i$, the set of 
sequences of the form $c_n = \sum_i \lambda_i^{\sigma^{-n}} z_i$ 
satisfies a twist-recurrence relation. Clearly we may alter the set 
as needed to ensure that the $z_i$ are linearly independent over
$R_0$, and that the $R_0$-module they span is saturated (that is, 
if $px$ lies in the span, so does $x$), or in other words,
the images of the $z_i$ in $R_0/\pi R_0$ are linearly independent over $\FF_p$.
If $\{z_1, \dots, z_k\}$ is the 
resulting set, then the equation
\[
\det \begin{pmatrix}
z_1 & z_1^\sigma & \cdots & z_1^{\sigma^{k}} \\
\vdots & & & \vdots \\
z_k & z_k^\sigma & \cdots & z_k^{\sigma^{k}} \\
c_n & c_{n+1}^\sigma & \cdots & c_{n+k}^{\sigma^{k}}
\end{pmatrix} = 0
\]
is satisfied when $c_n = \sum_i \lambda_i^{\sigma^{-n}} z_i$, and 
expanding this determinant along the bottom row gives us an equation 
of the form $d_0 c_n + d_1 c_{n+1}^\sigma + \cdots + d_k 
c_{n+k}^{\sigma^k} = 0$; all that remains to be checked is that $d_0$ 
and $d_k$ are not divisible by $\pi$.

Since $d_0 = (-1)^k d_k^\sigma$, we need
only check that $d_k$ is not divisible by $\pi$. In fact,
since $z_1, \dots, z_k$ are linearly independent mod $\pi$,
\[
d_k = \det \begin{pmatrix}
z_1 & z_1^\sigma & \cdots & z_1^{\sigma^{k-1}} \\
\vdots & & & \vdots \\
z_n & z_n^\sigma & \cdots & z_n^{\sigma^{k-1}} 
\end{pmatrix}
\equiv \det \begin{pmatrix}
z_1 & z_1^p & \cdots & z_1^{p^{k-1}} \\
\vdots & & & \vdots \\
z_n & z_n^p & \cdots & z_n^{p^{k-1}} 
\end{pmatrix}
\not\equiv 0 \pmod{\pi}.
\]
\end{proof}
Though we will not need this extra generalization, it is worth noting
that (by descent) the lemma also holds over
$W(L)$ for $L$ a perfect but not algebraically closed field.

We conclude with a remark about twist-recurrence in the case of 
primary interest in this paper, where $R = W_m(K)$ is the ring of 
Witt vectors of length $m$ over an algebraically closed field $K$ of 
characteristic $p$, with the canonical Frobenius.
It is natural to ask whether twist-recurrence of a 
sequence of Witt vectors is related to twist-recurrence of the 
sequence of $k$-th components for $k=1, \dots, m$; the following lemma 
answers this question affirmatively.
\begin{lemma} \label{split}
Let $\{c_0, c_1, \dots\}$ be a sequence of elements of $W_m(R)$, and
let $(w_{k,0}, \dots, w_{k,m-1})$ be the Witt vector corresponding to 
$c_k$.
\begin{enumerate}
\item
If $\{c_n\}$ is twist-recurrent, then so is $\{w_{n,i}\}$ for
$i=0,\dots, m-1$.
\item
If $\{w_{n,i}\}$ is twist-recurrent for $i=0, \dots, m-1$,
then so is $\{c_n\}$.
\end{enumerate}
Moreover, the derived twist-recurrence relations depend
only on the initial ones, and not on the particular sequences.
\end{lemma}
\begin{proof}
The key observation here is that a sequence $\{w_n\}$ of elements of
$K$ is twist-recurrent if and only if the sequence $\{[w_n]\}$ of
Teichm\"uller lifts of $\{w_n\}$ to $W_m(K)$ is twist-recurrent. 
Of course one implication is obvious.
For the other, note that in $W_m(K)$, $[x] =
t^{p^k}$ for any $k \geq m$ and any $t \in W_m(K)$ such that $t \equiv
x^{\sigma^{-k}} \pmod p$. Thus if we write $w_n = \sum_{i=1}^j
\lambda_i^{\sigma^{-n}} z_i$, we may select $y_i \in W_m(K)$ such that
$y_i \equiv z_i^{\sigma^{-m}} \pmod p$, and then
$\sum_i \lambda_i^{\sigma^{-n-m}} y_i$ is a lift of $w_n$. Now
\[
[w_n] = \left(\sum_i \lambda_i^{\sigma^{-n-m}} y_i\right)^{p^m}
= \sum_{e_1+\cdots+e_j=p^m} \frac{p^m!}{e_1!\cdots e_j!}
(\lambda_1^{e_1}\cdots
\lambda_j^{e_j})^{\sigma^{-n-m} }
y_1^{e_1}\cdots y_j^{e_j}
\]
satisfies a twist-recurrence relation not depending on the
$\lambda_i$.

Given the observation, the rest of the proof is straightforward. On
one hand, if the sequences of components of the Witt vectors are
twist-recurrent, so are the sequences of their Teichm\"uller lifts,
and the given sequence is simply a linear combination of these. On the
other hand, if the sequence of Witt vectors is twist-recurrent, then
the sequence of unit components is as well, as is its corresponding
sequence of Teichm\"uller lifts. Subtracting this sequence off and
dividing by $p$ gives a sequence of Witt vectors over $W_{m-1}(K)$
which is twist-recurrent, and the claim follows by induction.
\end{proof}

Beware that if $\{w_n\}$ is a twist-recurrent sequence in $K$,
the sequence $\{[w_n]\}$ of Teichm\"uller lifts to $W(K)$ may not be
twist-recurrent. 

\section{Completed algebraic closures}
Let $K$ be an algebraically closed field of characteristic $p$.
We use the $p$-adic ring $W(K)[[t^\QQ]]$ to explicitly construct the
$p$-adic completion of the algebraic closure of $W(K)$; in so doing, 
we will exploit the results of \cite{me}, though it is certainly 
possible to give independent derivations of the present results. 

We begin with the ring $W(K)[[t^\QQ]]$ and its $p$-adic completion, which may be described
as the set of series
$\sum_{i \in \QQ} c_i t^i$ such that for each $n \in \NN$, $\{i \in \QQ: v(c_i) \leq n\}$
is well-ordered. This ring is isomorphic to $W(K[[t^\QQ]])$;
we fix the isomorphism which sends $t \in W(K)[[t^\QQ]]$ to the
Teichm\"uller lift of $t \in K[[t^\QQ]]$, and use this isomorphism to identify the two rings.
Since $K[[t^\QQ]]$ is algebraically closed, it contains an algebraic closure of $K[[t]]$,
which we hereafter refer to as ``the'' algebraic closure $\overline{K[[t]]}$ of $K[[t]]$.
Likewise, since $W(K)[[p^\QQ]]$ is algebraically closed, it contains an algebraic closure
of $W(K)$, which we hereafter refer to as ``the'' algebraic closure $\overline{W(K)}$ of
$W(K)$.

There are a number of natural maps between the aforementioned rings,
which are summarized in the diagram below. Our first goal is to
establish the existence of the two dotted arrows.  (This relationship
between $\overline{K[[t]]}$ and $\overline{W(K)}$ is well-known; for
example, it occurs in the construction of the ``big rings'' of
Fontaine.)

\[
\xymatrix@=0pt@R=10pt@C=0pt@W=0pt@H=0pt{
     W(\overline{K[[t]]})\, \ar@{^{(}->}[rr]\ar@{-->>}[dd]\ar@{->>}[dr]
       & & W(K)[[t^\QQ]]^\wedge \ar@{->>}[dd]^(.35)\pi \ar@{->>}[dr] &
    \\ 
    & \overline{K[[t]]}/(t) \, \ar@{<-->}[dd] \ar@{^{(}->}[rr] &
     & K[[t^\QQ]]/(t) \ar@{<->}[dd]^(.35){\tilde{\pi}} 
    \\ 
     \overline{W(K)}^\wedge \, \ar@{^{(}->}[rr] \ar@{->>}[dr]
       & & W(K)[[p^\QQ]] \ar@{->>}[dr] &
    \\ 
     & \overline{W(K)}^\wedge/(p)\, \ar@{^{(}->}[rr] &
     & W(K)[[p^\QQ]]/(p)
    }
\]

\begin{theorem} \label{reduce}
The map $\pi: W(K)[[t^\QQ]]^\wedge \to W(K)[[p^\QQ]]$ induces a surjection of
$W(\overline{K[[t]]})$ onto $\overline{W(K)}^\wedge$ (and so
$\pibar$ induces an isomorphism of
$\overline{K[[t]]}/(t)$ with $\overline{W(K)}/(p)$).
\end{theorem}

The argument relies on a lemma which assert that the roots of a polynomial vary continuously
with the coefficients, in a manner independent of characteristic. The analogous statement for 
two polynomials over a single DVR, in which $k$ is not restricted, is well-known.

\begin{lemma} \label{lem}
Let $P(x) = x^n + a_1 x^{n-1} + \cdots + a_n$ be a polynomial over $K[[t^G]]$
and $Q(x) = x^n + b_1 x^{n-1} +
\cdots + b_n$ a polynomial over $W(K)[[p^G]]$
with the same Newton polygon.
Suppose $k \in (0, 1]$ has the property that for
for $i=1,\dots,n$, $\pibar(a_it^{-v_i}) \equiv b_i p^{-v_i} \pmod{p^{k}}$, where
$v_i$ is the $y$-coordinate of the point of the Newton polygon with $x$-coordinate
$i$.
Then for each root $y$ of $P$ of slope $s$, there exists a root $z$ of $Q$ of slope $s$
such that $\pibar(y) \equiv z \pmod{p^{k/m+s}}$ (where $m$ is the multiplicity of $s$
in the Newton polygon), and vice versa.
\end{lemma}
\begin{proof}
Let $s_1, \dots, s_l$ be the distinct slopes of the Newton polygon and $m_1, \dots, m_l$
their multiplicities. Then we can factor $P(x)$ as $P_1(x) \cdots P_l(x)$,
where $P_i(x)$ is the monic polynomial of degree $m_i$ whose roots are the roots
of $P$ which have valuation $s_i$ (counting multiplicities). Similarly we can factor $Q(x)$
as $Q_1(x) \cdots Q_l(x)$. Now observe that
\begin{align*}
P_l(x t^{-s_l}) &\equiv x^{m_l} + (a_{n-m_l+1}/a_{n-m_l}) x^{m_l-1} +
\cdots + (a_{n}/a_{n-m_l}) \pmod{t^k} \\ Q_l(x p^{-s_l}) &\equiv
x^{m_l} + (b_{n-m_l+1}/b_{n-m_l}) x^{m_l-1} + \cdots +
(b_{n}/b_{n-m_l}) \pmod{p^k}
\end{align*}
and so
\[
\pibar(P_l(x p^{-s_l})) \equiv Q_l(x p^{-s_l}) \pmod{p^k}.
\]
Moreover, $P/P_l$ and $Q/Q_l$ also obey the conditions of the
lemma. Thus by descending induction, we have $\pibar(P_i(x p^{-s_i}))
\equiv Q_i(x p^{-s_i}) \pmod{p^k}$ for $i=1, \dots, l$.  In other
words, the proof of the statement of the lemma reduces to the case in
which $P$ and $Q$ have only one slope, which we may assume is 0. In
this case, $n=m$ and $y_i=0$ for all $i$.

Let $y \in K[[t^G]]$ be a root of $P(x)$,
and choose $z \in W(K)[[p^G]]$
such that $\pibar(y) \equiv z \pmod{p}$. 
Let
\begin{align*}
c_i &= \sum_j a_{i-j} \binom{n-i+j}{j} y^j \\
d_i &= \sum_j b_{i-j} \binom{n-i+j}{j} z^j,
\end{align*}
so that $P(x-y) = \sum c_i x^{n-i}$ and $Q(x-z) = \sum d_i x^{n-i}$.
Then clearly $\pibar(c_i) \equiv d_i \pmod{p^{k}}$.  By the assumption
that $y$ is a root of $P(x)$, we have $c_n \equiv 0 \pmod{p^k}$.
That means that the Newton polygon of $Q(x-z)$ has at least one slope greater than
or equal to $k/m$, which is to say
$z$ is congruent to a root of $Q(x)$ modulo $p^k$, as desired.
The converse implication follows by an analogous argument.
\end{proof}

\begin{proof}[Proof of the Theorem]
%
%
%
We first establish that $\pi(W(\overline{K[[t]]}))$ is algebraically
closed; that is, for any $a_1, \dots, a_n \in W(\overline{K[[t]]})$,
the polynomial $R(x) = x^n + \pi(a_1) x^{n-1} + \cdots + \pi(a_n)$ has
a root in $\pi(W(\overline{K[[t]]}))$. Since this assertion only depends on
the values of $\pi(a_i)$, we may rechoose the $a_i$ to ensure that
\[
v_t(\rho(a_i)) = v_p(\pi(a_i)),
\]
where $\rho$ denotes reduction modulo $p$.
Then the polynomials
$P(x) = x^n + \rho(a_1) x^{n-1} + \cdots + \rho(a_n)$ over $\overline{K[[t]]}$
and $Q(x) = x^n + \pi(a_1) x^{n-1} + \cdots + \pi(a_n)$ over $W(K)[[p^\QQ]]$
clearly satisfy the conditions of Lemma~\ref{lem} with $k=1$. Thus for some root $z$
of $Q$ of highest slope $s$, there exists a root $y$ of $P$ in $K[[t^\QQ]]$,
also of slope $s$,
with $\pibar(y t^{-s}) \equiv z p^{-s} \pmod{p^{1/n}}$. Since the coefficients of $P$ lie
in the algebraically closed domain $\overline{K[[t]]}$, $y$ does as well.

Choose a lift $z$ of $y$ to $W(\overline{K[[t]]})$, making sure that
$z$ is divisible by $t^{v_t(y)}$.
We may reapply the above steps to $R(x-\pi(z))$ to get another ``approximate root'',
and so on. Each step increases the largest slope of the residual polynomial by
at least $1/n$, so the approximate roots converge in $W(K)[[p^\QQ]]$ to a root of
$R(x)$. Moreover, since we made $z$ divisible by $t^{v_t(y)}$, the approximate
roots also converge $(p,t)$-adically in $W(\overline{K[[t]]})$, which proves the claim.

In particular, we now have that $\pi(W(\overline{K[[t]]}))$ contains
$\overline{W(K)}^\wedge$, and it remains to establish that for any $r \in
W(\overline{K[[t]]})$, $\pi(r) \in \overline{W(K)}^\wedge$.
Without loss of generality, we may assume $v_p(\pi(r)) = v_t(\rho(r)) = 0$.

We will prove the claim by constructing a sequence $\{r_j\}$ of elements of $W(\overline{K[[t]]})$
and monic polynomials $P_j(x) = x^n + a_1^{(j)} x^{n-1} + \cdots + a_n^{(j)}$,
for $j=0, \dots, n$, such that for $j=0, \dots, n-1$:
\begin{enumerate}
\item[(a)] $r_0 = r$;
\item[(b)] $\pi(a_i^{(j)}) \in \overline{W(K)}^\wedge$; 
\item[(c)] $v_t(\rho(a_i^{(j)})) = v_p(\pi(a_i^{(j)}))$;
\item[(d)] $\rho(r_j)$ is a root of $\rho(P_j)$ of slope at least $j/n$;
\item[(e)] there exists a root $s_j$ of $\pi(P_j)$ such that $\pi(r_j) \equiv s_j \pmod{p^{(j+1)/n}}$;
\item[(f)] $\pi(r_{j} - r_{j+1}) = s_j$ and $P_{j+1}(x) = P_j(x+r_{j}-r_{j+1})$.
\end{enumerate}
Note that (b) can be arranged for $j=0$ by lifting a polynomial over $K[[t]]$
having $\rho(r)$ as a root (making sure that (c) is also satisfied). The fact that
(e) holds is a consequence of Lemma~\ref{lem}.

Observe that $\pi(r_0 - r_n) \in \overline{W(K)}^\wedge$ and $v_p(r_n) \geq 1$.
That is, $\pi(r)$ is congruent to an element of $\overline{W(K)}^\wedge$ modulo
$p$. Repeating this process with $r_n$ in place of $r$ exhibits succesive
elements of $\overline{W(K)}^\wedge$ congruent to $\pi(r)$ modulo $p^i$ for
all $i \in \NN$, which implies that $\pi(r) \in \overline{W(K)}^\wedge$. Thus
$\pi(W(\overline{K[[t]]})) = \overline{W(K)}^\wedge$, as desired.
\end{proof}

This gives a rather nonexplicit construction of the completed algebraic closure of 
$W(K)$. We explicate the construction 
in two steps, first giving a concrete description of
$W(\overline{K[[t]]})$.
In order to effect this description, we recall some definitions from 
\cite{me}, adapted slightly to accommodate our new definition of 
twist-recurrence. For $a, b$ in $\NN$, we define the set $S_{a,b} 
\subseteq \QQ$ as follows:
\[
S_{a, b} = \left\{ \frac{1}{a}(n - b_1 p^{-1} - b_2 p^{-2} - \cdots):
n \in \NN, b_i \in \{0, \dots, p-1\}, \sum b_i \leq b \right\}.
\]
We also define $T_b = S_{1,b} \cap (0, 1)$.

For $R$ a $p$-adically complete ring with Frobenius, we say
a function $\map{f}{T_b}{R}$ is
\emph{twist-recurrent} 
if there exists $k, l \in \NN$ and $d_0, \dots, d_{k-1} \in R$, with 
$d_0$ not divisible by $p$, such that the twist-recurrence relation
(\ref{twist}) holds
for any sequence $\{c_n\}$ of the form
\begin{equation} \label{seq}
c_n = f(- b_1 p^{-1} - \cdots - b_{j-1} p^{-j+1} -
p^{-n+l}(b_j p^{-j} + \cdots) ) \quad (n \geq 0)
\end{equation}
for $j \in \NN$ and $b_1, b_2, \dots \in \{0, \dots, p-1\}$ with $\sum b_i
\leq b$. (The explicit insertion of $l$ is required because eventually
twist-recurrent sequences need not be twist-recurrent.)

\begin{theorem}
The ring $W_n(\overline{K[[t]]})$
is isomorphic to the set of $x \in W_n(K)[[t^\QQ]]$ such that
the following conditions hold.
\begin{enumerate}
\item
There exist $a, b \in \NN$ such that
if $x_i \not\equiv 0 \pmod{p^n}$, then $i \in S_{a,b}$.
\item
For each nonnegative integer $m$, the function $\map{f_m}{T_b}{W_n(K)}$ given by
$f_m(i) = x_{(m+i)/a}$ is twist-recurrent.
\item
The functions $f_m$ span a $W_n(K)$-module of finite rank (in the 
space of all maps from $T_b$ to $W_n(K)$ with the natural $W_n(K)$-module 
structure).
\end{enumerate}
\end{theorem}
\begin{proof}
It suffices to note that the latter set 
is a ring (by Lemma~\ref{twistcomb}), and has residue ring $\overline{K[[t]]}$
(by \cite{me}).
\end{proof}

With this in hand, we now give a concrete description of the completed 
algebraic closure of $W(K)$.
Let $R = W(K)[[p^\QQ]] = W(K)[[t^\QQ]]/(t-p)$ be the $p$-adic ring 
over $W(K)$ with value group $\QQ$.
As noted in Section~\ref{sec:def}, each element $x$ of $R$
admits a unique representation in $W(K)[[t^G]]$ of the form $\sum
[x_i] t^i$, where $a_i \in K$ and brackets denote the Teichm\"uller
map. Beware that the $x_i$ do not exhibit quite the same behavior as
their counterparts in $K[[t^G]]$; most notably, they are not additive.

Let $B$ denote the subset of $R$ consisting of those $x$ satisfying
the following conditions for each $n \in \NN$:
\begin{enumerate}
\item
There exist $a, b \in \NN$ such that $x$ is
supported on $S_{a,b} \cup (n, \infty)$.
\item
For some (any) choice of $a,b$ as above and for each $j < n$, the 
function $\map{f}{T_b}{K}$ given by (\ref{seq}) is twist-recurrent.
\end{enumerate}

As in \cite{me}, a
simplification occurs in the case $k = \Fpbar$:
the twist-recurrent condition can be
replaced by the simpler condition that for some $M, N \in \NN$, the
sequences become periodic after $M$ terms with period length at most $N$.

\begin{theorem}
The image of $W(\overline{K[[t]]})$ under the quotient $W(K)[[t^\QQ]] 
\to W(K)[[p^\QQ]]$ is equal to $B$ (which thus is
the $p$-adic completion of an algebraic
closure of $W(K)$).
\end{theorem}
\begin{proof}
The most difficult part of the proof is showing that $B$ is actually a ring!
If $z=x+y$ then
\[
z_i = x_i + y_i + \sum_{j = 1}^\infty
P_j(x_i,y_i, x_{i-1}, y_{i-1}, \dots, x_{i-k}, y_{i-k})
\]
for 
certain universal polynomials $P_j$.
Suppose that $x$ and $y$ are twist-recurrent.
To check $z \in B$, we must verify the condition given above for each $n$.
However, for given $n$, we may replace the infinite sum over $j$ by a sum
with $j \leq n$ without changing the values of $z_i$ for $i<n$. This
expresses $z_i$ as a polynomial in the $x_{i-k}$ and $y_{i-k}$, and
any fixed polynomial in sequences satisfying given twist-recurrence
relations satisfies one as well (by Lemma~\ref{twistcomb}). Thus $z$ is
twist-recurrent.

If $z = xy$, then
\[
z_i = \sum_{j+k=i} \sum_{m=0}^\infty Q_m(x_j,y_k,x_{j-1},y_{k-1},\cdots, x_{j-m}, y_{k-m})
\]
for certain universal polynomials $Q_m$.
The argument that $z$ is twist-recurrent given that $x$ and $y$ are is similar
to that for multiplication, except one must additionally 
note that as over a sequence of the form
(\ref{seq}), the number of terms in the outer sum is \emph{uniformly} bounded. 
See the analogous section of \cite{me} for a detailed description of why this is the case.

Given that $B$ is a ring, which evidently is $p$-adically complete,
let $R$ be the image of $W(\overline{K[[t]]})$ in $W(K)[[p^\QQ]]$.
To show $R=B$, we need only show $R/pR = B/pB$ as
as subsets of $W(K)[[p^\QQ]]/(p) = K[[t^\QQ]]/(t)$.
But this follows from the description of $\overline{K[[t]]}$ given in \cite{me}, together
with Lemma~\ref{split}.
\end{proof}

Among other things, this result resolves a question of Lampert 
\cite{lamp}, who
asked for the possible order types of the support of an element of the $p$-adic
ring $\ZZ_p[[p^\QQ]]$
algebraic over $\ZZ_p$. The theorem shows that, writing $\omega$ 
for the first countable ordinal, the support of an algebraic element 
has order type at most $\omega^\omega$ (whereas the support of a
general element of $\ZZ_p[[p^\QQ]]$ may have any countable order type).
This observation may be useful
for implementing the extraction of roots of polynomials over $\ZZ_p$
on a computer.

It is worth noting explicitly that the mixed-characteristic result 
above is weaker than the equal-characteristic analogue, in 
that we deduce no criterion for identifying $\overline{W(K)}$ within 
its completion. Indeed, such a criterion would 
have to be more subtle in this case:
for example, in \cite{me} it is shown that every
algebraic series over $K[[t]]$ is supported on $S_{a,b}$ for some
$a,b$, but this already appears to fail over $\ZZ_p$ for the polynomial $x^p -
p^{p-1}x = p^{p-1}$. This question clearly deserves to be better 
understood, but we will not comment further on it here.

It is also worth reiterating an observation from \cite{me}: a 
twist-recurrent series can be specified modulo $p^n$ by a finite number of 
elements of $K$. This may make twist-recurrent series amenable to machine
computations.

\end{document}